\documentclass[12pt,reqno]{amsart}

\newcommand{\sz}{{s_0}}

\newcommand{\wt}{\widetilde}
\newcommand{\wh}{\widehat}
\newcommand{\sm}{\setminus}
\newcommand{\PP}{{\mathbb P}}
\newcommand{\R}{{\mathbb R}}
\newcommand{\C}{{\mathbb C}}

\newcommand{\Z}{{\mathbb Z}}

\renewcommand{\d}{\partial}

\newcommand{\codim}{{\operatorname{codim\,}}}
\newcommand{\supp}{{\operatorname{Supp\,}}}
\newcommand{\rank}{{\operatorname{rank}}}
\newcommand{\loc}{{\operatorname{loc}}}
\renewcommand{\phi}{\varphi}
\newcommand{\al}{\alpha}
\newcommand{\be}{\beta}

\newcommand{\Ga}{\Gamma}

\newcommand{\ep}{\varepsilon}
\newcommand{\de}{\delta}
\newcommand{\De}{\Delta}

\newcommand{\Om}{\Omega}

\newcommand{\bcal}{\mathcal{B}}
\newcommand{\ccal}{\mathcal{C}}

\newcommand{\gcal}{\mathcal{G}}
\newcommand{\jcal}{\mathcal{J}}
\newcommand{\lcal}{\mathcal{L}}
\newcommand{\ocal}{\mathcal{O}} 
\newcommand{\pcal}{\mathcal{P}}
\newcommand{\qcal}{\mathcal{Q}}

\newtheorem{theo}{{\sc Theorem}} 
\newtheorem{cor}[theo]{{\sc Corollary}}
\newtheorem{lem}[theo]{{\sc Lemma}}
\newtheorem{prop}[theo]{{\sc Proposition}}

\newenvironment{rem}{\medskip\noindent{\it Remark:\/} }{\medskip}

\title[Compact singularities of meromorphic mappings]{Compact
singularities of meromorphic mappings between complex 3-dimensional
manifolds}
\author{Sergei Ivashkovich}
\address{U.F.R. de Mathematiques, Universit\'e de Lille-I, 59655 Villeneuve 
d'Ascq Cedex, France.}
\email{ivachkov@gat.univ-lille1.fr}
\author{Bernard Shiffman}
\address{Department of Mathematics, Johns Hopkins University, Baltimore,
MD
21218, USA} 
\email{shiffman@math.jhu.edu}

\thanks{Research of the second author partially supported by NSF grant
\#DMS-9800479.}

\date{\today}

\begin{document}

\begin{abstract}We prove that a meromorphic map defined on the complement of a 
compact subset of a three-dimensional Stein manifold $M$ and with values in a
compact complex three-fold $X$ extends to the complement of a finite set of
points.  If $X$ is simply connected, then the map extends to all of $M$.
\end{abstract}

\maketitle

\section*{Introduction}

The study of the extendibility of holomorphic and meromorphic mappings began
with the classical theorem of Hartogs \cite{Ha} (see \cite{Si}).
\begin{quote}
{\it Let $K$ be a compact subset of a domain $M \subset \C^n$, $n\ge 2$,
such that $M 
\setminus K$ is connected, and let $f:M \setminus K\to \C$ be a
holomorphic function. Then there exists a holomorphic function $\hat
f:M \to \C $  extending $f$, i.e., $\hat f\mid_{M \setminus K}=f$.
}\end{quote}
Shortly after Hartogs proved his theorem, E. E. Levi \cite{Le} discovered that
this extension result holds true also for meromorphic functions.  

A natural problem is to understand under what conditions Hartogs' Theorem
(respectively Levi's Theorem) holds when the mapping $f$ takes values in  a
general complex manifold
$X$ rather than $\C$ (respectively
$\C\PP^1$). Of course it is immediate that Hartogs' Theorem 
remains valid for holomorphic mappings with values in a Stein manifold $X$,
since such a manifold $X$ can be embedded into
$\C^N$. It similarly follows that Levi's Theorem also remains 
valid for meromorphic mappings with values in compact projective manifolds.

In 1971, Griffiths \cite{Gr} and the
second author
\cite{Sh} independently showed that Hartogs' Theorem is
valid for holomorphic mappings into manifolds $X$ carrying a complete
Hermitian metric  with non-positive holomorphic scalar curvature, answering a
question was asked by  Chern in \cite{Che}. Concerning the meromorphic
mapping problem, the first author \cite{Iv-4} proved that Hartogs extension
holds for meromorphic maps into  compact K\"ahler manifolds.

We  recall  two more results here due to K. Stein 
and M. Chazal. Stein proved in \cite{St} that Hartogs' Theorem holds for holomorphic maps if $\dim
X\le n-2$. Recently Chazal \cite{Ch} relaxed this condition to $\dim X\le 
n-1$ and more generally $f$ can be meromorphic.  The next case of interest is
the equidimensional case $\dim X = n$.  It is well known that one doesn't
always have meromorphic extension in this case, as is illustrated by the
(holomorphic) projection $f:\C^n\to X=\C^n/\Z$ to the Hopf  manifold. (The
$\Z$-action  is given by $z\buildrel {n}\over\mapsto 2^nz$.)  The goal of this
paper is to  show that, at least for dimension $\le 3$, the singularity at 0
of the Hopf map $f$ is the only type of singularity that can occur for
equidimensional meromorphic maps:

\begin{theo}\label{maintheorem} Let $K$ be a compact set with connected
complement in a Stein manifold $M$ of dimension $3$, let $X$ be a
compact complex manifold of the same dimension and let
$f:M\sm K\to X$ be a meromorphic map.  Then there exists a finite
set $\{a_1,\ldots ,a_d\}\subset K$ such that $f$ has a meromorphic
extension
$\hat f:M \sm \{a_1,\ldots ,a_d\}\to X$, and if $ B(a_j)$ are
disjoint coordinate balls centered at $a_j$, then $\hat
f(\d  B(a_j))$ is not homologous to zero in $X$ ($1\le j\le d$).
\end{theo}

The same result when both $M $ and $X$ have dimension two  was 
proved by the first author in \cite[Corollary~4(b)]{Iv-2}.  It is open whether
this result if valid for equidimensional maps of dimension greater than 3. Of
course, one cannot expect to obtain such results when the dimension of
$X$ is greater  then the dimension of $M$; see the remark in
\S 1 below.

In the case of the Hopf map $f:\C^3\setminus \{0\}\to \C^3/\Z$ mentioned above,
$A=\{0\}$.  Of course, the elements
$\hat f(\d B(a_j))$ of the fifth homology group are very special; they are
often called spherical shells.  If, for example,
$\hat f$ is a holomorphic embedding in a neighborhood of $\d B(a_j)$, then
$X$ is of a very restricted type:  it is a deformation of the Hopf
3-fold (see \cite{Ka-1}).\medskip

In particular, if the singular set $A$ is nonempty, then
$H^5(X,\R)\neq 0$.   
Poincar\'e duality then implies the following:

\begin{cor} If under the conditions of Theorem \ref{maintheorem},
$H^1(X,\R)=0$, then $f$ extends meromorphically to all of
$M$.\end{cor}

\section {Reductions}

For degenerate mappings,
the result is known and is due to F. Chazal \cite{Ch}. Hence, in the
sequel, we always suppose that $f$ is nondegenerate; i.e.,
${\rm rank}\,f=3$. 

We let $ \De(r)=\{z\in\C: |z|<r\}$ denote the disk or radius $r$
about $0$, and we write $\De=\De(r)$.  We consider the polydisk
$\De^n(r)=\De(r)^n$ and 
``annulus" $A^n(r,1)=\De^n \sm \overline{\De^n(r)}$. We 
shall make frequent use of the following Hartogs figure in
$\C^3$:
$$H^2_1(r)= \big[\Delta(1-r)\times \De^2\big] \cup\big[ \Delta \times
A^2(r,1)
\big]\,.$$

By the standard method of extending analytic objects
(see for example
\cite{Si}), it suffices to prove either of the following two equivalent
results:

\begin{prop}\label{main} Let $X$ be a compact complex $3$-fold and let
$f:H^2_1(r)\rightarrow X$ be a nondegenerate meromorphic map. Then
there is  a discrete set $\{a_j\}\subset \Delta^3
\setminus H^2_1(r)$  and a meromorphic extension $\hat f:\Delta^3
\sm \{a_j\}\to X$ such that if  $ B(a_j)$ are
disjoint balls in $\Delta^3$ centered at $a_j$, then $\hat
f(\d  B(a_j))$ is not homologous to zero in $X$.\end{prop}

\begin{prop}\label{psconcave} Let $M,\ W$ be open sets in
$\C^3$, and suppose $p\in M\cap\d W$ such that $W$ is smooth and 
strictly pseudoconvex at $p$. Let $U=M\sm \overline W$. Suppose $f:U\to X$
is a nondegenerate meromorphic map to a compact $3$-fold $X$.  Then
there is an open set
$\wt U\supset U\cup\{p\}$ such that
$f$ has a meromorphic extension $\hat f$ to
$\wt U\sm \{p\}$, and  either $\hat f$ is meromorphic at $p$ or
$\hat
f(\d  B(p))$ is not homologous to zero in $X$, where $ B(p)$
is a ball about $p$ contained in $\wt U$.\end{prop}

\smallskip We note that Proposition~\ref{psconcave} follows from Proposition
\ref{main} with the additional simplifying assumption that $f$ is  {\it
holomorphic} on 
$\De \times A^2(r,1)$.  To
see this, we  observe that the set of points of indeterminacy
$I_f$ of our meromorphic  map $f$ has codimension at least two, i.e., is a
curve together with a  discrete set of points. Let
$M,\,W,\,p$ be as in Proposition~\ref{psconcave} and
let $\De^2_{p_1}=\{p_1\}\times\De^2$ denote the vertical bidisk passing through
$p=(p_1,p_2,p_3)$.  We can assume, after making a quadratic change of
coordinates, that
$\De^2_{p_1}\cap I_f$ contains no curves and $\De^2_{p_1}\cap
\overline W=\{p\}$.  After translating and
stretching coordinates, we then obtain a
Hartogs figure
$H_1^2(r)$  contained in $U$, with $p$ in the corresponding polydisk
$\De^3$, such that $[\De \times A^2(r,1)]\cap I_f=\emptyset$, so we can
apply the modified Proposition~\ref{main} to obtain the conclusion of
Proposition~\ref{psconcave} with $\wt U=\De^3$.  

Proposition \ref{main} follows from Proposition~\ref{psconcave}, since the 
Hartogs figure can be exhausted by a family of strictly pseudoconcave 
hypersurfaces and this family can be continued to exhaust $\Delta^3$. Thus,
when proving Proposition \ref{main}, we may assume that $f$ is  holomorphic on 
$\De \times A^2(r,1)$.

\smallskip
\begin{rem} The reader may observe that these results involves 
extension from a ``$1$-concave" 3-dimensional domain. It is  worthwhile
to note that in general there is no extension of meromorphic maps with 
values in compact 3-dimensional manifolds from $2$-concave domains, such
as the classical Hartogs figure 
  
$$H^1_2(r):= \big[\Delta^2(1-r)\times \Delta \big] \cup\big[ \Delta^2
\times A^1(r)
\big]\,.$$
Namely, M. Kato constructed in \cite{Ka-2} an example of a compact
complex three-fold 
$X$ and a holomorphic mapping $f:\C^2\setminus  \bar B\to X$ defined on
the complement of a ball $B\subset \C^2$, such that every point of the
sphere $\d B$ is an essential singular point of
$f$. 
\end{rem}

\medskip We shall prove Proposition \ref{main} in \S 2 after
we make the following reductions: 

\begin{enumerate}
\item[a)] First of all, as was already explained, we can assume that 
$f:H_1^2(r)\to X$ is nondegenerate and holomorphic on 
$\De \times A^2(r,1)$.

\item[b)] We can further assume that there is no hypersurface in $H^2_1(r)$
which $f$ sends to a point.  If such a hypersurface exists, then by shrinking
$H^2_1(r)$ a little bit, we can suppose that there are finitely many of them.
Then by blowing up the image points sufficiently many times, we obtain a
modification $\wh X$ of $X$ together with a lift $\hat f$ of $f$ to a
meromorphic map $\hat f:U\rightarrow \wh X$ having the desired property.
After extending $\hat f$, we can push it down to extend $f$ itself.

\item[c)] We write
$$A^2_s=\{s\}\times A^2 (r, 1)\,,\quad s\in\De\,.$$
After again shrinking $H^2_1(r)$ a little, we can suppose  that $A^2_s$
contains no curves contracted by
$f$ to a point, for all $s\in\De$.  Indeed, since ${\rm
rank}\,f=3$, there are at most 1-parameter families of contracted curves. 
We consider small quadratic changes of the $z_1$ coordinate:
$\wt z_1 = z_1 + Q(z_2,z_3)$, where $Q$ is a polynomial of degree 2 with
small coefficients.   The set of
such Q such that
$\wt z_1$ is constant on a fixed holomorphic curve is of codimension at
least 2. Whence, for an open
dense set of such
$Q$, the coordinate function
$\wt z_1$ is nonconstant on each contracted curve; i.e., for all
$s\in\De$,
$\wt z_1^{-1}(s)$ contains no contracted curves.

\item[d)] By the above argument, we can also assume that none of the
$A^2_s$ are contained in the critical set $C_f$ of $f$.

\item[e)] By the argument below, we can also assume
that for all $s\in\De$, there do not exist nonempty disjoint open
subsets $V_1,V_2$ of $A^2_s$ with $f(V_1)=f(V_2)$.

\end{enumerate}

To show that we can realize property (e) after a change of coordinates, we
let
$$U=\big[\De\times A^2(r,1)\big]\sm (I_f \cup C_f)$$ and we consider the
set
$$D=
\{(z,w)\in U\times U: z\ne w,\ f(z)=f(w)\}\,,$$ which is an analytic
subvariety of $U\times U$ minus the diagonal.  Note that $D$ is locally
given as the graph of a biholomorphic map (and thus is a smooth
3-dimensional submanifold).
It suffices to show that we can make a small perturbation of
coordinates so that \begin{equation}\label{dimD} \dim D\cap
(\De^2_s\times \De^2_s) \le 1\end{equation} for all $s\in\De$.

To show (\ref{dimD}), we let
$\pcal^l_n$ denote the vector space of polynomials of degree $\le l$ on
$\C^n$. Note that 
\begin{equation}\label{dim} \dim
\pcal^l_n=\left(\!\!\begin{array}{c}l+n\\
n\end{array}\!\!\right)\,.\end{equation}  We also let
$\jcal ^l_a(g)\in \pcal^l_n$ denote the
$l$-jet of a germ
$g\in\,{}_n\ocal_a$ ($a\in\C^n$).
We shall use the following lemma:

\begin{lem}\label{jets} Let $\phi:\De^m\to \C^n$ be a holomorphic map such
that $\phi(0)=a$ and $\rank_0\phi=m$. Then 
$$\codim_{\pcal^l_n}\{f\in \pcal^l_n :\jcal^l_0(f\circ \phi)=0\}=
\left(\!\!\begin{array}{c}l+m\\
m\end{array}\!\!\right)\,.$$
\end{lem}

\begin{proof} By a change of coordinates, we can assume without loss of
generality that \break $\phi(z_1,\dots,z_m)=(z_1,\dots,z_m,0,\dots,0)$. The
result then follows immediately from (\ref{dimD}).\end{proof}

We let $\qcal^l$ denote the set of polynomials $g$ in $\pcal^l_3$ such that
$dg$ does not vanish on $\De^3(2)$.  (Note that small polynomial
perturbations of the coordinate function $z_1$ are in $\qcal^l$.)  Let
$a=(z_0,w_0)\in D\cap\De^6(2)$ be arbitrary, and let $\bcal_a$
denote the set of polynomials $g$ in $\qcal^5$ with $$\dim_a\{(z,w)\in
D:g(z)=g(z_0), g(w)=g(w_0)\}>1\,.$$ We shall show that
\begin{equation}\label{codimB} \codim_{\qcal^5}\bcal_a\ge
6\,.\end{equation}   Since
$\dim D=3$, (\ref{codimB}) implies that we can choose  $g\in\pcal^5_3$
such that $g$ is a small deformation of the coordinate function $z_1$ and 
$$ \dim D\cap \left(g^{-1}(s)\times g^{-1}(t)\right) \le 1\quad \mbox{for
all\ } s,t\in\De\,.$$ If we then replace $z_1$ with $\tilde z_1=g$,
(\ref{dimD}) will be satisfied.

To verify (\ref{codimB}), we first consider an arbitrary quadratic
polynomial
$g_1\in\qcal^2$, and we let
$$E=\{(z,w)\in D: g_1(z)=g_1(z_0)\}\,.$$ Since $D$ is locally given as a
graph and 
$dg_1(z_0)\ne 0$,  $E$ is smooth at $z_0$.  Now let
$\phi=(\phi_1,\phi_2):\De^2\to E$ with
$\phi(0)=a$ and $\rank_0\phi=2$.  This implies that
$\rank_0\phi_1=\rank_0\phi_2=2$.  Let
$\bcal'(g_1)$ denote the set of
$g_2\in\qcal^2$ such that $\jcal^2_0 ( g_2 \circ \phi_2)=0$. By Lemma
\ref{jets}, $\codim \bcal'(g_1)\ge {4\choose 2}=6$.  Furthermore, we note
that if we replace $g_1$ with a germ $\tilde g_1\in {}_3\ocal_{z_0}$ with
the same 2-jet at $z_0$, then we can choose $\wt \phi_2$ with  the same
2-jet (at 0) as $\phi_2$ so that 
$(\phi_1,\wt
\phi_2):\De^2\to E$ has the same 2-jet (at 0) as $\phi$. Thus if
$g_2\in\qcal^2\sm\bcal'(g_1)$, we have
$\jcal^2_0 ( g_2
\circ
\wt\phi_2)=\jcal^2_0 ( g_2
\circ \phi_2)\ne 0$.  Furthermore, if we also replace this $g_2$ with 
$\tilde g_2\in {}_3\ocal_{w_0}$ with the same 2-jet at $w_0$, then
$\jcal^2_0 ( \tilde g_2
\circ \wt\phi_2)\ne 0$, and hence $$\dim_a\{(z,w)\in D:\tilde
g_1(z)=g_1(z_0),\;\tilde g_2(w)=g_2(w_0)\}=1\,.$$ 

Now consider the linear map
$$\tau_a:\pcal^5_3\to\pcal^2_3\times
\pcal^2_3\,,\quad g\mapsto \left(\jcal^2_{z_0}(g),
\jcal^2_{w_0}(g)\right)\,.$$ By the above discussion,
$\bcal_a
\subset
\tau^{-1}_a(\bcal'_a)$, where
$$\bcal'_a=\{(g_1,g_2):g_1\in\qcal^2,\;g_2\in\bcal'(g_1)\}\,.$$ Since
$\tau$ is surjective, it follows that
$$\codim_{\qcal^5} \bcal_a\ge\codim_{\qcal^2\times\qcal^2}\bcal'_a \ge
6\,.$$
This completes the verification of
(\ref{dimD}), and hence conditions (a)--(e) above can be satisfied.

\bigskip

\section{Proof of Proposition \ref{main}}

We are now prepared to prove the Hartogs extension property.  By our
construction above, we may assume that the map
$f:H^2_1(r)\to X$  of Proposition
\ref{main} possesses the following properties:

\begin{enumerate}
\item[i)] $f$ is non-degenerate and holomorphic on a neighborhood of
$\overline{\Delta \times A^2(r,1)}$;

\item[ii)] for all $s\in \Delta$, the set
$A^2_s$ contains no curves contracted by
$f$ to a point;

\item[iii)] for all $s\in\De$, there do not exist nonempty disjoint open
subsets $V_1,V_2$ of $A^2_s$ with $f(V_1)=f(V_2)$;

\end{enumerate}

We must show that $f$ extends meromorphically to $\De^3$ minus a
discrete set of points. Denote by $W$ some open subset of $\Delta$ such
that $f$ can be meromorphically extended onto the Hartogs
domain $$H_W(r):=\big[W\times
\De^2\big]
\cup \big[
\Delta \times A^2 (r,1)\big]\,.$$

Let $\Omega$ be a strictly positive $(2,2)$-form on $X$ with
$dd^c\Omega=0$. Existence of such a form on the compact 3-dimensional 
manifold $X$ follows from the absence of nonconstant plurisubharmonic 
functions on $X$ via duality and the Hahn-Banach theorem. In fact even more
is  true. Every Hermitian metric on $X$ is conformally equivalent to a
metric whose  associated $(1,1)$-form $\omega $ is $dd^c$-closed, see
\cite{Ga}. In the  sequel, we shall take $\Omega =\omega^2$, where $\omega
$ is such a Gauduchon form.  Denote by 
$T$ the pull-back of $\Omega$ by
$f$, i.e.\
$T=f^*\Omega$. More accurately, $f^*\Omega $ is defined in the case of
meromorphic
$f$ as follows. Let $\wt\Ga_f$ denote the desingularization of
the graph $\Ga_f\subset H_W(r)\times X$ of $f$ and let $\pi_1:\wt\Ga_f\to
H_W(r)$ and 
$\pi_2:\wt\Ga_f\to X$  be the natural projections. Note that 
$\pi_1$ is proper by the very definition of 
meromorphic map. Define 
\begin{equation}\label{T}f^*\Omega:=\pi_{1*}\pi_2^*\Omega\,.\end{equation}
The current $T=f^*\Omega$ is a positive bidegree $(2,2)$ current on
$H_W(r)$. Being the push-forward of a smooth form (on a
desingularization of $\Gamma_f$), $T$  has coefficients in
$\lcal^1_\loc(H_W(r))$.

 To see that the push-forward of a smooth form
$\eta$ by a modification $\pi:\wt M\to M$ has coefficients in $\lcal^1_\loc$,
it suffices
to show that $\pi_*\eta$ has no mass on the center $C$ of $\pi$.  (In our
case $C=I_f$.) But for any test form
$\phi$ on $M$ and any sequence $\rho_n\to \chi_C$ with
$0\le\rho_n\le 1$, we have $(\pi_*\eta, \rho_n \phi) = \int_{\wt M}
\eta\wedge \pi^*(\rho_n \phi)\to 0$. Hence $\|\pi_*\eta\|(E)=0$. (In fact,
this holds when $\pi$ is any surjective holomorphic map that is proper on
$\supp\eta$.)

It follows immediately from (\ref{T}) that $dd^cT=0$. Moreover,
$T$ is smooth on 
$H_W(r)\setminus  I_f$, since outside the set $I_f$ of indeterminacy
points of $f$, it is the usual pull-back of the smooth form $\Om$.
   
We write $\De^2_s:=\{s\}\times \De^2$, for $s\in\De$.
The function
$$
\mu(s) := \int_{\De^2_{s}}f^* \Omega$$
is well defined for all $s\in W$, since by the
above, $(f^* \Omega)|_{\De^2_s}=(f_{\De^2_s})_*\Omega$ is a positive,
bidegree $(2,2)$-current on a neighborhood of $\overline{\De^2_s}$ and is
in $\lcal^1_\loc$. We remark that
$\mu(s)$ is nothing but the  volume of
$f(\De^2_{s})$ with respect to $\Omega $ counted with multiplicities.

\begin{lem}  The function $\mu$ is positive and  smooth on $W$, and its 
Laplacian
$\Delta\mu$ smoothly extends onto the whole unit disk $\Delta $.\end{lem}

\begin{proof} (We follow the method of proof of \cite[Lemma~3.1]{Iv-1}.)
The positivity of
$\mu$ follows from the positivity of
$f^*
\Omega$ and property (ii) above. To show that $\De\mu$ extends to the unit
disk, we begin by writing
$$T =  \sum_{\al,\be =1}^3 t_{\alpha
\bar\beta} dz^{\stackrel{\alpha}{\vee}} \wedge d\bar
z^{\stackrel{\beta}{\vee}}\,,$$ where
$dz^{\stackrel{1}{\vee}}=dz_2\wedge dz_3,\ 
dz^{\stackrel{2}{\vee}}=-dz_1\wedge dz_3,\ 
dz^{\stackrel{3}{\vee}}=dz_1\wedge dz_2$. The function $\mu$ is then given
by  $$\mu(z_1)  = 
\int_{\De^2} t_{1\bar1} (z_1, z_2, z_3) dz_2\wedge d\bar z_2\wedge
dz_3\wedge d\bar z_3\,.$$
Let
$$
T^{\varepsilon }=  \sum_{\al,\be =1}^3
t_{\alpha
\bar\beta}^\ep dz^{\stackrel{\alpha}{\vee}} \wedge d\bar
z^{\stackrel{\beta}{\vee}}
$$
be the smoothing of $T$ by convolution; the $T^{\varepsilon }$ are smooth 
forms converging to $T$ in $\lcal^1$ as $\varepsilon \to 0$. On
$H_W(r)\setminus  I_f$ the convergence is in the $\ccal^{\infty }$ topology.
The functions
$$\mu^\ep(z_1): =\int_{\De^2_{z_1}} T^\ep =
\int_{\De^2_{z_1}} t_{1\bar1}^\ep (z_1, z_2, z_3) dz_2\wedge d\bar
z_2\wedge dz_3\wedge d\bar z_3$$
are smooth in $W$. 

The condition $dd^cT=0$ reads as 

$$
\sum_{\alpha ,\beta}
{\d^2t_{\alpha \bar\beta }\over \d z_{\alpha }\d \bar z_{\beta }}=0.
$$
So,\begin{eqnarray*}
\Delta\mu^\ep\, (z_1)&=&  4
\int_{\De^2_{z_1}} {\d^2 t_{1\bar1}^\ep\over \d z_1\d \bar z_1}  dz_2\wedge
dz_3\wedge d\bar z_2\wedge d\bar z_3\\ &=& 
-4 \int_{\De^2_{z_1}}
\sum_{(\alpha,\beta )
\not=(1,1)} {\d^2t_{\alpha \bar\beta }^\ep\over \d z_{\alpha }\d \bar
z_{\beta }}
 dz_2\wedge d\bar
z_2\wedge dz_3\wedge d\bar z_3\\ &=& -4
\int_{\d \De^2_{z_1}} \left[  \sum_{\alpha =2}^3 \pm
{\d t_{\alpha \bar1}^\ep\over \d \bar z_1} dz_{5-\alpha}\wedge d\bar
z_2
\wedge d\bar z_3 + \sum_{\alpha=1}^3\sum_{\beta 
=2}^3 \pm
{\d t_{\alpha \bar\beta }^\ep\over \d z_{\alpha }}dz_2\wedge dz_3\wedge
d\bar z_{5-\beta}\right]\,.
\end{eqnarray*}
Since $f$ is holomorphic on a neighborhood of $\De\times\d \De^2$, the current
$T$ is smooth on $\De\times\d \De^2$ and thus
$\Delta\mu^{\varepsilon }$ converges
smoothly  to the function $\psi$ given by
\begin{equation}\label{psi}\psi(z_1)= -4
\int_{\d \De^2_{z_1}} \left[  \sum_{\alpha =2}^3 \pm
{\d t_{\alpha \bar 1}\over \d \bar z_1} dz_{5-\alpha}\wedge d\bar
z_2
\wedge d\bar z_3 + \sum_{\alpha=1}^3\sum_{\beta 
=2}^3 \pm
{\d t_{\alpha \bar\beta }\over \d z_{\alpha }}dz_2\wedge dz_3\wedge
d\bar z_{5-\beta}\right]\,.
\end{equation}
But (\ref{psi}) defines a smooth function on all of $\De$.
While $\mu^{\varepsilon }\to
\mu $ in 
$\lcal^1$ on $W$, so $ \Delta\mu^{\varepsilon }\to \Delta \mu $ on $W$. This
shows  that 
$\Delta\mu =\psi$ is smooth and smoothly extends onto the disk $\Delta $.  Thus
$\mu $ is also smooth on $W$.
\end{proof}

\begin{lem} Suppose that $f$ is non-degenerate and that there exists a
sequence $\{s_n\}\in W$ converging to $s_0\in\De$ such that $\mu (s_n)$ is
bounded.  Then:

{\rm 1)}  $f_0:=f|{A^{2}_{s_0}}$ meromorphically extends onto
$ \De^2_{s_{0}}$;

{\rm 2)} the volumes of the graphs $\Gamma_{f_{s_{n}}}$ are uniformly bounded
in $n$;

{\rm 3)}  $f$ meromorphically extends onto $U_0 \times  \De^2$ for some
neighborhood $U_0$ of $s_0$.
\label{mubound}\end{lem}

\begin{proof} 
1) We let $f_n= f|\De^2_{s_n}$ and we write $F_n=f_n( \De^2_{s_n})$. We
further write 
$\Sigma_n=f_n(\d  \De^2_{s_n}),\;$ $\Sigma_0=f_0(\d  \De^2_{s_0})$. Since
the volumes
$\mu(s_n)$ of the $F_n$  are uniformly bounded,
by Bishop's theorem (see for example
\cite{HS}) we can assume, after passing to a subsequence, that $F_n$
converges to a pure 2-dimensional analytic subset $F$ of $X\sm \Sigma_0$. Note
that $\bar F=F\cup \Sigma_0$.

\medskip\noindent {\it Case 1.\/} $\bar F$ is a subvariety of $X$.

In this case
$$f_0|_{A^2_{s_0}}: A^2_{s_0}\to F^0$$
is a holomorphic map into an irreducible component $F^0$ of $\bar F$. If
$\wh F^0$ denotes a desingularization of $F^0$, then $f_0$ lifts to a
meromorphic map $\hat f_0$ from $A^2_\sz$ to $\wh F^0$.  By the 2-dimensional
version of Theorem \ref{maintheorem} proved in
\cite[Cor.~4(b)]{Iv-2}, $\hat f_0$ extends meromorphically onto $ \De^2_\sz$
minus a finite set of points. If this set is nonempty, then $\hat
f_0(\d \De^2_\sz)$ would not be homologous to zero  in $\wh F^0$, and hence
$\Sigma_0=f_0(\d \De^2_\sz)$ would not be homologous to zero in $F^0$.  But 
$\Sigma_0=\lim \Sigma_n=\lim \d F_n =\d \lim F_n$ in the sense of currents.
Since $\supp \lim F_n\subset \bar F$,
$\Sigma_0$ is homologous to zero in $\bar F$ and therefore in $F^0$, 
a contradiction.  So $f_0$ extends onto all of
$ \De^2_\sz$.

\medskip\noindent {\it Case 2.\/}  $\bar F$ is not a subvariety of $X$.

Let $F^0$ be the irreducible component of $F$ containing $f_0(A^2_\sz)
\setminus \Sigma_0$.  Define the analytic space $E=F_0\cup A^2_\sz/\!\!\sim$,
where the equivalence relation is defined as follows: The points $a\in F^0$
and $b\in A^2_\sz$ are equivalent iff $a=f_0(b)$, and necessarily $b'\in
A^2_\sz$ is equivalent to $b$ iff $f_0(b)=f_0(b')$. By property (ii) above,
this is a proper equivalence relation and hence $E$ is a complex space. Let
$\pi:E\to \bar F^0 \subset X$ be the projection defined by $\pi(a)=a$ for
$a\in F^0$ and $\pi(b)=f_0(b)$ for $b\in A^2_\sz$.  Let $\wh E{\buildrel
{\eta}\over \to} E$ denote the normalization.  By property (iii), the map
$f_0:A^2_\sz\to E$ is generically one-to-one and thus is a normalization of
its image.  By the uniqueness of the normalization, $f_0$ lifts to a map $\hat
f_0:A^2_\sz\to\wh E$, i.e., $\eta\circ \hat f_0 =f_0$.  The map $\hat f_0$
is a biholomorphism onto its image.

The boundary $\d\wh E$, being biholomorphic to $\partial \De^2$, is strictly
pseudoconvex after shrinking slightly, so by Grauert's theorem, $\wh E$ can
be blown down to a normal Stein space.  This easily yields an extension of
$f_0$ onto $ \De^2_\sz$.

2) We denote the extension of $f_0$ onto $ \De^2_\sz$ also by $f_0$.  
Let $F'$ be the maximal compact pure 2-dimensional variety contained in
$\Bar F$.  (In Case 1 above, $F'=\bar F$, whereas in Case 2,
$F'=\overline{F\setminus F^0}$.) 
We consider the pure two-dimensional analytic set
$$
\Gamma = \Gamma_{f_0}\cup (I_{f_0}\times F')\,.
$$
in $(\Delta\times \De^2)\times X$, where $I_{f_0}$ is the (finite) set of
points of indeterminacy of $f_0$.

\smallskip\noindent {\sl Step 1.} {\it We claim that for all $\ep >0$ the graph 
$\Gamma_{f_n}$ belongs to the
$\varepsilon $-neighborhood of
$\Gamma$, for $n \gg 0$.
}

\smallskip
  Neighborhoods are taken with respect to the 
Euclidean metric on $\C^3$ and Gauduchon metric on $X$. (In fact,
any choice of metric works as well as this one.)  This claim follows
immediately from Lemma \ref{sequence} below.

We say that a sequence of meromorphic maps $f_n:U\to X$ converges to a
holomorphic map $f_0$ on a domain $U$ if for all compact subsets $K\subset
U$, $I_{f_n}\cap K=\emptyset$ for $n
\gg 0$ and $f_n \to f_0$ uniformly on $K$.

\begin{lem} \label{sequence} Let $f_n:\bar \De^2 \to X$ be a sequence of
meromorphic maps, where
$X$ is a compact complex manifold.  Suppose that $f_n$ is holomorphic on
$A$, where
$A=A^2(r,1)$. If there exists a meromorphic map
$f_0:\De^2 \to X$ such that $f_n|_A \to
f_0|_A$, then
$f_n \to f_0$ on $\De^2 \setminus I_{f_0}$. 
\end{lem}

\smallskip\noindent Lemma \ref{sequence} is a special case of Proposition
1.1.1  in \cite{Iv-3}. (Proposition
1.1.1  in \cite{Iv-3} is stated in terms of ``strong convergence" of
meromorphic maps.  However, if $\{f_n\}$ strongly converges to a holomorphic
map, then the sequence converges in the above sense.  This is the content of
the ``Rouch\'e principle" of  \cite[Theorem~1]{Iv-3}.)

\bigskip
To complete the proof of (2), we consider a point $p\in \Gamma $ and take any
open $W\ni p$ adapted to $\Gamma $, i.e.  biholomorphic to  $\Delta^2\times
\Delta^4=U\times B$ in such a way  that $(\bar U\times \d B)\cap \Gamma
=\emptyset $. Then for $n\gg 1$, we have 
$\Gamma_{f_n}\cap  (\bar U\times \d B)=\emptyset $ and thus
$p\mid_{\Gamma_{f_n}}:   
\Gamma_{f_n}\cap (U\times  B)\to U$ is a $d_n$-sheeted analytic 
covering, where $p:U\times B\to U$ is a natural projection.

\smallskip\noindent {\sl Step 2.} {\it  The number  $\{ d_n\} $ of sheets is 
uniformly bounded.
}
\smallskip

Consider the following two cases.
\smallskip\noindent {\sl Case 1.} $p\in (I_{f_0}\times F')\setminus
\Gamma_{f_0}.$ In this case as $W\ni p$ we can take the following
neighborhood. Let $p=(a,b)$, where $a\in I_{f_0}\subset \C^3$ and $b\in F'
\subset X$. Take a neighborhood of $b$ in $X$ of the form $\Delta^2\times
\Delta $ such that $F'\cap (\bar \Delta^2\times \d\Delta )=\emptyset $. Then
take some small $\Delta^3\ni a$ in $\C^3$ and put $U=\Delta^2$ and
$B=\Delta\times \Delta^3$. If the number $d_n$ of sheets of the analytic cover
$\pi_U:\Gamma_{f_n}\cap (U\times B)\to U$ is not bounded, it will contradict
the fact that $f_n(\De^2_{s_n})\cap (\Delta^2\times\Delta )=F_n\cap
(\Delta^2\times\Delta )$ has uniformly bounded volume (counted with
multiplicities).

One should remark now that boundedness of the number of sheets does not
depend on the particular choice of the adapted neighborhood of $p$.

\smallskip\noindent {\sl Case 2.} $p\in \Gamma_{f_0}$. Let $W=U\times B\ni p$ 
be some adapted neighborhood. Find a point $q\in U$ such that all its
pre-images 
$\{ q_1,...,q_N\} = \pi_U^{-1}(q)\cap \Gamma $ are smooth points of $\Gamma $ and 
$\pi_U$ is a biholomorphism between neighborhoods $V_j\ni q_j$ on $\Gamma $ 
and $V$ on $U$. Denote by $b_j$ the projection of $q_j$ into $B$. Take 
mutually disjoint polydisks $B_j\subset B$  with centers $b_j$. Consider 
$W_j:=V_j\times B_j$ as adapted neighborhoods of $\Gamma $ in $q_j$. They 
are adapted also for $\Gamma_{f_n}, n>>0.$ Denote by $d^j_n$ the 
corresponding
number of sheets. If $d_n$ is not bounded then at least one sequence $d^j_n$ 
is also unbounded. Fix $j$ with $d^n_j$ unbounded. 

If $q_j\in (I_{f_0}\times F')\setminus \Gamma_{f_0}$, then
everything  reduces to Case 1. So let $q_j\in \Gamma_{f_0}$.
Perturbing $q$ and thus $q_j$ if necessary,  we can suppose that $q_j$ is
a point where our map 
$f$ is holomorphic. More precisely $q_j=(a,f(a))$ for some $a\in
\Delta\times \De^2 
\subset \C^3$. Now the contradiction is immediate, because the graphs 
$\Gamma_{f_n}$ uniformly approach $\Gamma_{f_0}$ while $f_n$
converges  to $f$ in a neighborhood of $a$.

\medskip\noindent
3) We are exactly under the assumptions of Proposition 1.3 of
\cite{Iv-2}, i.e., we can apply the ``Continuity Principle." (The condition of 
boundedness of the cycle geometry is insured by Proposition 1.4 from
\cite{Iv-2}.) This gives us an extension of
$f$ onto
$U_\sz
\times  \De^2$.
\end{proof}

\bigskip
Let us proceed further with the proof of the theorem. Let $W$ be the maximal
open  subset of the disc $\De$ such that $f$ meromorphically extends onto
$H_W(r)$. 

\begin{lem} $\De\setminus W$ is a closed complete polar set in $\De$.
\end{lem}

\medskip
The proof is the same as that of Lemma 2.4 from \cite{Iv-2} and will be
omitted. 

\bigskip It suffices to show that there exists
$\hat f:\De^3(1-\de)\to X$ satisfying the conclusion of Proposition~\ref{main}
for arbitrary $\de>0$. We now repeat the above arguments using  two slightly
deformed coordinate systems
$(z'_1 ,z_2,z_3)$ and $(z''_1 ,z_2,z_3)$, where
$$z'_1=z_1+\ep z_2 +O(|z|^2)\,,\quad z''_1=z_1+\ep z_3 +O(|z|^2)
\,.$$
Here the $O(|z|^2)$ terms are chosen so that conditions (i)--(iii) at the
beginning of this section are satisfied for each of the two coordinate
systems, after shrinking $r$ if necessary.  (As was shown in \S 1, these terms
can be taken to be polynomials consisting of terms of degrees 2 through 5.) We
choose
$\ep$ and the
$O(|z|^2)$ terms to be small enough so that $\De^3(1-\de)\subset \De'{}^3
\subset \De^3,\;$ $\De^3(1-\de)\subset \De''{}^3 \subset \De^3$, where
$\De'{}^3$ and $\De''{}^3$ are the polydisks of radius
$1-\frac{\de}{2}$ in the new coordinates.   

Applying the above argument to the new coordinate systems, we obtain maximal
open $W',\; W''$ in $\De\tilde{}:=\De(1-\frac{\de}{2})$ such that $f$ extends
meromorphically to the Hartogs domains $H'_{W'}(r),\; H''_{W''}(r)$.  We let
$S_1=\De\tilde{}\setminus W,\;S_2=\De\tilde{}\setminus W',\;
S_3=\De\tilde{}\setminus W''$.  Now consider the coordinates 
$$w_1=z_1\,,\quad w_2=z'_1\,,\quad w_3=z''_1$$ and let $U$ denote the image of
$\De^3(1-\de)$ under the coordinate map $(w_1,w_2,w_3)$ . (We may assume that
$z'_1,z''_1$ are chosen so that the $w_j$ indeed provide coordinates on $\De^3$.)

In terms of the $w$-coordinates, $f$ then extends to a meromorphic map $\hat
f$ on
$U\setminus  (S_1\times S_2\times S_3)$.  Now let $s_0$ be an arbitrary point
in $S:=S_1\times S_2\times S_3$.  We must show that $s_0$ is an isolated point
of
$S$ and that $\hat
f(\d  B_{s_0}(r))$ is not homologous to zero in $X$, for any ball $B_{s_0}(r)$
centered at $s_0$ such that 
$ B_{s_0}(r)\cap S=\{s_0\} $.  

Since polar sets in $\C$ are
of  Hausdorff dimension zero, we can choose a polydisk $\Delta^3_0 $ about
$s_0$ such  that the set
$K:=S\cap \De^3_0$ is compact. An identical proof to that of Lemma 3.3 from
\cite{Iv-1} now shows that the  current $T=f^*\Omega $ has
locally  summable coefficients on all of $\De^3_0$. Hence $T$  
extends to a unique current $\tilde T$ on $\De^3_0$ with $\lcal^1_{\rm loc}$
coefficients.  The following lemma then tells us that $dd^c\tilde T$ is of
order 0:

\smallskip
\begin{lem} {\rm \cite[Proposition 2.3]{Iv-1}} Let $K$ be a complete
pluripolar, compact set in a strictly  pseudoconvex domain $D \subset \C^n$
and $T$ a positive, bidimension 
$(1,1)$ current in $D \setminus K$. Suppose that:

\smallskip\noindent
{\rm 1)} $dd^cT\le 0$ in $D \setminus K$,

\noindent
{\rm 2)} $T$ has locally finite mass in a neighborhood of $K$,

\noindent
{\rm 3)} $dT$ and $d^cT$ have measure coefficients on $D
\setminus K$.

\smallskip\noindent Then the current $dd^c\tilde T$ has measure
coefficients in
$D$.
\end{lem}
\noindent(Condition (3) on  $dT$ and $d^cT$ was
omitted in \cite{Iv-1}, but is used  in the proof.  In our case
$T=\pi_{1*}\pi^*_2\Om$, so this condition  follows from the fact
that
$dT=\pi_{1*}\pi^*_2d\Om $  is  the push-forward of a smooth form
by a proper map, and similarly for $d^cT$.)

Since $dd^cT=0$, the support of the current $dd^c\tilde T$ must be contained
in $K$.  We also conclude from the  Lemma 2.6 in \cite{Iv-2}
that  $dd^c\tilde T\le 0$. Thus we can write
$dd^c\tilde T= \nu 
\omega_e^3$, where $\omega_e$ is the Euclidean  K\"ahler form on $\C^3$.
Then for any ball $B_{s_0}(r)\subset\!\subset \De^3_0$ about $s_0$ with $\d
B_{s_0}(r)\cap K=\emptyset$, we have  that either 

\begin{equation}
\nu (K\cap B_{s_0}(r))=0, 
\end{equation}
or 
\begin{equation}\label{negative}\begin{array}{lcl}
0&>&\displaystyle\nu (K\cap B_{s_0}(r))\ =\ \int_{B_\sz (r)}dd^c\tilde T\ =\ 
\lim_{\varepsilon \to 
0}\int_{B_{s_0}(r)}dd^c\tilde T_{\varepsilon }\\[14pt]
&=&\displaystyle \lim_{\varepsilon \to 0}\int_{\partial B_\sz (r)}d^c\tilde
T_{\varepsilon }\ =\ 
\int_{\partial B_{s_0}(r)}d^cT_{\varepsilon }\ =\ \int_{f(\partial
B_{s_0}(r))} d^c\Omega\, .\end{array}\end{equation}

\smallskip\noindent\sl
Case 1. $\nu (K\cap B_{s_0}(r))=0.$

\smallskip\rm In this case the negativity of $\tilde T$ implies that 
$dd^c\tilde T=0$. Therefore we can find a polydisk 
neighborhood $\De^3\subset B_{s_0}(r)$ of $s_0$ and a $(2,1)$-form $\Gamma $
in  $\De^3$ such that:

\smallskip\noindent
1) $f$ is holomorphic in a neighborhood of $\De \times \partial \De^2$ (and 
therefore $\tilde T$ is smooth there);

\noindent
2) $\tilde T=i(\partial \bar\Gamma -\bar\partial \Ga)$ in a neighborhood of 
$\bar \De^3$;

\noindent
3) $\Gamma $ is smooth in a neighborhood of  $\De \times \partial \De^2$.

\smallskip The zero-dimensionality of $K$ implies that there exists a nonempty 
open $W\subset \De $ such that $f$ is defined and meromorphic on $W\times
\De^2$ and that $s_0^1\in \partial W\cap \De $, where $s_0^1$ is the first
coordinate  of $s_0$. As before we let
$$
\mu (z_1)=\int_{\De^2_{z_1}}\tilde T=i\int_{\partial \De^2_{z_1}}(\bar\Gamma -
\Gamma )\,. 
$$
By the smoothness of $\Gamma $, the function $\mu$ is bounded.
Therefore by Lemma \ref{mubound},
$f$ extends meromorphically to a neighborhood of $s_0$.

\smallskip\noindent\sl
Case 2. $\nu (S\cap B_{s_0}(r))<0.$

\smallskip\rm  
By (\ref{negative}), the 5-cycle $f(\partial B_{s_0}(r))$ is not homologous to
zero in
$X$. Furthermore, $\int_{ f(\partial B_{s_0}(r))}d^c\Omega$ depends only 
on the integer homology class of 
$f(\partial B_{s_0}(r))$, since $dd^c\Om=0$. Hence, 
$$\int_{ f(\partial B_{s_0}(r))}d^c\Omega \le - \de<0\,,$$
where
$\de$ is independent of $s_0$ and $r$ (and depends only on $X$ and $\Om$).
This shows that $K$ is finite,
and completes the proof.
\qed

\smallskip Remark that our proof gives more. Namely, if $\Sigma \subset M
\setminus \{ a_1,...,a_d\} $ is not homologous to zero in $ M \setminus \{
a_1,...,a_d\} $ then $f(\Sigma )$ is not homologous to zero in $X$.

\medskip
\section {Generalizations and open questions}

\medskip In \cite{Iv-2}, the classes ${\pcal}_k^-$ and ${\gcal}_k$ of complex
spaces were introduced. Recall that ${\pcal}_k^-$ is the class of normal
complex spaces which carry a strictly positive $(k,k)$-form $\Omega^{k,k}$
with $dd^c\Omega^{k,k}\le 0$, and ${\gcal}_k$ is the subclass of ${\pcal}_k^-$
which consists of complex spaces carrying a strictly positive $(k,k)$-form
$\Omega^{k,k}$ with $dd^c\Omega^{k,k}=0$.  Note that ${\gcal}_k$ contains all
compact complex manifolds of dimension $k+1$.

It is easy to observe that our above proof gives the following more general
statement of Proposition~\ref{main}:

\begin{prop} Let $X$ be a compact complex manifold in the class ${\pcal}_2^-$. 
Then every meromorphic map $f:H_1^2(r)\to X$ extends meromorphically onto
$\Delta^3\setminus A$, where $A$ is a closed, complete pluripolar subset of
Hausdorff dimension zero. If moreover $X\in {\gcal}_2$ then $A$ is discrete
and for every ball $B$ with center $a\in A$ such that $\d B\cap A=\emptyset $,
$f(\d B)$ is not homologous to zero in $X$.
\end{prop}

To consider the extension of mappings from higher dimensional domains, we
introduce the Hartogs figures
$$H^k_d(r):= \big[\Delta^d(1-r)\times \Delta^k \big] \cup\big[ \Delta^d \times
A^k(r) \big]\subset \C^{d+k}\,.$$ We conjecture that the analogous result
should hold for meromorphic mappings from $H_d^k(r)$ to compact manifolds (and
spaces) in the classes ${\pcal}_k^-$ and ${\gcal}_k$.  In particular,
Theorem~\ref{maintheorem} should be true for meromorphic mappings between
equidimensional manifolds in all dimensions. The main difficulty lies in the
fact that it is impossible in general to make the reductions (a)--(c) of \S
1. (Note that reductions (d)--(e) can be achieved in all dimensions.)
However, these reductions are unnecessary in the case when our map is locally
biholomorphic, as we state below.

\begin{prop} Let $X$ be a compact complex space of dimension $k+1$.
Then every holomorphic map $f:H_1^k(r)\to X$ with zero-dimensional fibers
extends meromorphically onto $\Delta^{k+1}\setminus A$, where $A$ is discrete,
and for every ball $B$ with center $a\in A$ such that $\d B\cap A=\emptyset $,
$f(\d B)$ is not homologous to zero in $X$.
\end{prop}
 
The proof is by induction on the dimension $n=k+1$. For the inductive step,
the function $\mu$ is defined in terms of the push-forward of a $dd^c$-closed,
positive
$(k,k)$-form $\Omega$ on a desingularization of $X$.

\end{document}